\documentclass{amsart}
\usepackage[russian]{babel}
\usepackage[cp1251]{inputenc}
\usepackage{geometry}
\usepackage{amsthm}

\newtheorem {theorem} {Теорема}

\newtheorem {claim} [theorem] {Утверждение}

\theoremstyle {definition}

\theoremstyle {remark}

\def\deg{\mathop{\fam0 deg}}
\def\Z{\Bbb Z}

\begin{document}

\title{О теореме Понтрягина-Стинрода-Ву}
\author{Душан Реповш, Михаил Скопенков и Фулвиа Спаггиари}
\address{Institute for Mathematics, Physics and Mechanics, University of
Ljubljana, P. O. Box 2964, 1001 Ljubljana,
Slovenia. E-mail: dusan.repovs@uni-lj.si.}
\address{Москва 119992, Московский Государственный Университет,
механико-математический факультет, кафедра дифференциальной
геометрии и приложений. E-mail: skopenkov@rambler.ru.}
\address{Dipartimento di Matematica, Universita` di Modena e
Reggio Emilia, via Campi 213/B, 41100 Modena, Italy. E-mail: spaggiari.fulvia@unimo.it.}

\thanks{Реповш частично поддержан Министерством Образования и Науки Республики Словения, программа No.~101-509.}
\thanks{Скопенков частично поддержан грантом ИНТАС 06-1000014-6277,
грантами Российского Фонда Фундаментальных Исследований 05-01-00993-a,
06-01-72551-НЦНИЛ-а, 07-01-00648-a, Грантом Президента Российской Федерации для государственной поддержки ведущих научных школ Российской Федерации,
проект НШ-4578.2006.1, программой Министерства Образования и Науки ''Развитие научного потенциала высшей школы'',
проект РНП 2.1.1.7988, Фондом поддержки молодых ученых ''Конкурс Мёбиуса''.}

%
%


\begin{abstract}
Данная работа посвящена гомотопической классификации
отображений $(n+1)$-мерных многообразий в $n$-мерную сферу. Для
отображения $f:M^{n+1}\to S^n$ {\it степень} $\deg f\in H_1(M;\Bbb Z)$ есть класс,
двойственный  $f^*[S^n]$, где $[S^n]\in H^n(S^n;\Z)$
--- фундаментальный класс. В работе приводится простое
доказательство следующего частного случая теоремы
Понтрягина-Стинрода-Ву:

\smallskip

{\bf Теорема.} {\it Пусть $M$ --- связное ориентируемое замкнутое
гладкое $(n+1)$-многообразие, $n\ge3$. Тогда отображение
$
\deg:\pi^n(M)\to H_1(M;\Bbb Z)
$
является

1-1 отображением (биекцией), если
$w_2(M)\cdot\rho_2 H_2(M;\Bbb Z)\ne0$;

2-1 отображением (то есть каждый элемент
$\alpha\in H_1(M;\Bbb Z)$ имеет ровно 2 прообраза)--- иначе.}

\smallskip

Доказательство основано на конструкции Понтрягина-Тома и
геометрическом построении классов Штиффеля-Уитни $w_i(M)$.
\end{abstract}

\maketitle

Пусть $M$ --- связное ориентируемое замкнутое гладкое многообразие
размерности $n+k$. Обозначим через $L_k(M)$ множество $k$-мерных
оснащенных зацеплений в $M$ с точностью до оснащенного кобордизма.
Согласно конструкции Понтрягина-Тома, множество $L_k(M)$ находится
в 1-1 соответствии с множеством  $\pi^n(M)=[M;S^n]$ непрерывных
отображений $M\to S^n$ с точностью до гомотопии. Основная цель
нашей статьи --- описать множество $L_k(M)=\pi^n(M)$ при $k=1$ в
''стабильном ранге'' $n\ge3$. В статьях \cite{Pon38, Ste47} (сравни с \cite[\S30.3]{FoFu89})
 описание множества $L_k(M)$ было сведено к вычислению
стинродовых квадратов, которое было проделано Ву (сравни с \cite[\S30.2.D]{FoFu89}).
В данной статье приводится простое доказательство этой
классификационной теоремы Понтрягина-Стинрода-Ву. Есть основания
полагать, что это и есть первоначальное рассуждение Понтрягина,
которое он не опубликовал, обратившись сразу к общему случаю,
когда $M$ является произвольным полиэдром (см. теорему~\ref{d-th2} ниже и
замечанием после ее формулировки).

Данная классификация оснащенных зацеплений основана на понятии
естественной ориентации и степени оснащеного зацепления,
определяемых следующим образом. Фиксируем ориентацию $M$.
Рассмотрим точку $x$ оснащенного зацепления $L$ и оснащение
$f_1,\dots,f_n$ в этой точке. Базис $e_1,\dots,e_k$ касательного
пространства $T_x(L)$ назовем {\it положительным}, если базис
$e_1,\dots,e_k,f_1,\dots,f_n$ пространства $T_x(M)$ положителен.
{\it Степень} $\deg L$ оснащенного зацепления $L$ есть
(целочисленный) гомологический класс цикла $L$ с естественной
ориентацией. Тем самым имеем отображение
$$
\deg:L_k(M)\to H_k(M;\Z).
$$
Теорема Хопфа-Уитни 1932--35 утверждает, что данное отображение
биективно при $k=0$ и сюръективно при $k=1$.

\begin{theorem}\label{d-th1}
(a) Пусть $M$ --- связное ориентируемое
замкнутое гладкое $(n+1)$-многообразие, $n\ge3$. Тогда
$\deg:L_1(M)\to H_1(M;\Z)$ является

1-1 отображением (то есть биекцией), если
$w_2(M)\cdot\rho_2 H_2(M;\Z)\ne0$;

2-1 отображением (то есть каждый элемент
$\alpha\in H_1(M;\Bbb Z)$ имеет ровно 2 прообраза) --- иначе.

(b) Пусть $M$ --- связное ориентируемое замкнутое гладкое
$(n+2)$-многообразие, $n\ge3$. Элемент $\alpha$ лежит в образе
отображения $\deg:L_2(M)\to H_2(M;\Z)$ если и только если
$\rho_2\alpha\cdot w_2(M)=0$.
\end{theorem}


Здесь $\cdot$ обозначает умножение $H_k(M;\Z_2)\times
H^k(M;\Z_2)\to\Z_2$. Но в доказательстве этой теоремы удобнее
заменить классы Штиффеля-Уитни двойственными к ним гомологическими
классами (их геометрическое определение напоминается ниже). Ниже
эти гомологические классы обозначаются теми же буквами $w_i$ и
$\bar w_i$, а $\cdot$ обозначает пересечение $H_i(M)\times
H_j(M)\to H_{i+j-m}(M)$. Интересно было бы узнать ответ на вопрос, при каких $n$ условие $\rho_2\beta\cdot
w_2(M)\ne0$  в случае (а) теоремы~\ref{d-th1} можно заменить на
$w_2(M)\ne0$; ср. с \cite{FU91}. 



\begin{theorem}[Понтрягин]\label{d-th2}
(a) Пусть $M^3$ --- связное
ориентируемое замкнутое гладкое 3-многообразие. Тогда для каждого элемента $\alpha\in
H_1 (M^3; \mathbb{Z})$ имеется взаимно-однозначное соответствие между множествами
$\deg\nolimits ^ {-1} \alpha$ и $\mathbb{Z}_{2d(\alpha)}$,
где $d (\alpha) $ --- делимость проекции элемента
$\alpha$ на свободную часть группы $H_1 (M^3; \mathbb{Z}) $.

(b) Пусть $M$ --- связное ориентируемое замкнутое 4-многообразие.
Элемент $\alpha$ лежит в образе отображения $\deg:L_2(M)\to
H_2(M;\Z)$ если и только если $\alpha\cdot\alpha=0$.

\end{theorem}

Теорема Понтрягина~\ref{d-th2}.b может быть доказана аналогично доказательству теоремы~\ref{d-th1}.b ниже.
Методы данной работы могут быть применены также для доказательства теоремы~\ref{d-th2}.a,
сформулированной без доказательства в \cite{Pon41}. (На самом деле, теорема~\ref{d-th2}.a была
сформулирована не в статье \cite{Pon41}, написанной по-английски, а только в резюме,
написанном по-русски, без указаний к доказательству.) Доказательство теоремы~\ref{d-th2} опубликовано,
например, в~\cite{CRS3}.

\smallskip



{\bf Геометрическое определение гомологических классов Штиффеля-Уитни.} 

Рассмотрим систему $s$ касательных векторных полей общего положения на $M$. Пусть
$\Sigma\subset M$ --- множество точек, в которых данные векторные поля линейно зависимы.
По трансверсальности \cite[\S10.3]{DNF79} $\Sigma$ является псевдомногообразием 
в $M$. {\it Класс
Штиффеля-Уитни} $w_{n+2-s}(L)\in H_{s-1}(M;\Z_2)$ --- это гомологический класс
псевдомногообразия 
$\Sigma$ (это первое препятствие к существованию системы $s$ линейно
независимых касательных векторных полей на $M$).

Данное определение легко обобщается на случай, когда касательные векторные поля на $M$
заменяются на сечения произвольного векторного расслоения с базой $M$. Для
подмногообразия $L\subset M$ классы Штиффеля-Уитни нормального расслоения к $L$ в $M$ и
ограничения $TM$ на $L$ определяются аналогично и обозначаются через $\bar w_i(L)$ и
$w_i(M)|_L$, соответственно.

Мы собираемся использовать также относительные классы Штиффеля-Уитни. Они определяются
следующим образом. Пусть $L\subset M\times I$ --- $l$-мерный (неоснащенный) кобордизм
между оснащенными подмногообразиями $L_1$ и $L_2$. Продолжим оснащение $f$ края $\partial
L=L_1\cup L_2$ до системы $m-l-1$ нормальных векторных полей общего положения на $L$.
Определим относительный класс Штиффеля-Уитни $\bar w_2(L,f)\in H_{l-2}(L;\Z_2)$ как класс
$(l-2)$-псевдомногообразия, 
в точках которого построенные векторные поля линейно зависимы
(это первое препятствие к продолжению оснащения $\partial L$ на $L$). Мы будем опускать
$f$ из обозначения классов Штиффеля-Уитни в случае, когда это не приведет к
неоднозначному прочтению.

\smallskip


{\it Доказательство теоремы~\ref{d-th1}.b.} (Д. Реповш и Ф. Спаггиари) Возьмем произвольный элемент $\alpha\in H_2(M;\Z)$. Реализуем
класс $\alpha$ ориентируемым 2-многообразием $L$. Очевидно, $\alpha\in\operatorname{Im}
\deg$ если и только если подмногообразие $L$ можно оснастить (при некотором выборе представителя $L$).

В дальнейшем считаем $L$ связным. Действительно, пусть некоторое несвязное $L$ реализует
класс $\alpha$ и $L$ можно оснастить. Тогда сделаем связную сумму всех компонент $L$.
Полученное подмногообразие связно, его также можно оснастить и оно реализует тот же
гомологический класс.

Покажем, что {\it $L$ можно оснастить, если только если $\bar w_2(L)=0$}. Данное условие,
очевидно, является необходимым, согласно данному выше геометрическому определению классов
Штиффеля-Уитни. Докажем достаточность этого условия. Предположим, что $\bar w_2(L)=0$.
Стандартными методами теории препятствий можно показать, что тогда на $L$ существует
ортонормированная система нормальных векторных полей $f_1,\dots,f_{n-1}$ (так как $n\ge3$
и $\dim L=2$ \cite{FoFu89}). Так как $L^2$ и $M^{n+2}$ ориентируемы, то нормальное расслоение
к $L$ в $M$ ориентируемо. Фиксируем его ориентацию. Дополняя $f_1,\dots,f_{n-1}$
единичным вектором $f_n$ до положительного базиса (относительно выбранной ориентации),
получим искомое оснащение многообразия $L$. Выделенное курсивом
утверждение доказано.

Теперь теорема~\ref{d-th1}.b следует из равенства
$$\bar w_2(L)=w_2(M)|_L=w_2(M)\cdot [L]=w_2(M)\cdot\rho_2\alpha.$$
Дадим необходимые пояснения. Здесь и далее мы отождествляем с $\Z_2$ все группы $H_0(X;\Z_2)$,
изоморфные $\Z_2$.
В данной цепочке равенств первое равенство получается из
формулы Ву для классов Штиффеля-Уитни суммы двух расслоений:
$w_2(M)|_L=w_2(L)+w_1(L)\cdot\bar w_1(L)+\bar w_2(L)$, в которой $w_2(L)=w_1(L)=0$, так
как $L$ --- ориентируемое 2-многообразие. (Первое равенство может быть также доказано
непосредственно.) Второе равенство следует из Геометрического определения выше, поскольку
$L$ связно.
Последнее равенство в нашей цепочке выполняется по определению подмногообразия $L$.
$\square$
\smallskip


{\it Начало доказательства теоремы~\ref{d-th1}.}
Рассмотрим элемент $\alpha\in H_1(M;\Z)$. Наша
цель --- выяснить, сколько существует попарно некобордантных оснащенных зацеплений $L_i$
таких, что $\deg L_i=\alpha$.

Пусть $L_1,L_2\subset M$ --- пара оснащенных зацеплений (1-подмногообразий), таких что
$\deg L_1=\deg L_2=\alpha$. Обозначим через $[L_1]$ и $[L_2]$ их классы в $L_1(M)$. Так
как $L_1$ и $L_2$ гомологичны, по общему положению существует вложенный ориентированный (неоснащенный)
кобордизм между ними
$$
L\subset M\times I,\partial L=L_1\sqcup L_2.
$$
Ясно, что $[L_1]=[L_2]$, если и только если для некоторого такого $L$ заданное оснащение
$\partial L$ продолжается на $L$. Так как $M$ связно, мы можем считать $L$ связным.

Покажем, что {\it  оснащение края $\partial L$ продолжается на все $L$ если и только если
относительный класс Штиффеля-Уитни $\bar w_2(L)=0$.}
Согласно геометрическому определению выше это условие
необходимо. Докажем его достаточность. Предположим, что $\bar w_2(L)=0$. Стандартными
методами теории препятствий можно показать, что тогда первые $n-1$ векторных полей
оснащения продолжаются до ортонормированной системы векторных полей на $L$ (так как
$n\ge3$ и $\dim L=2$ \cite{FoFu89}).
Так как $L$ и $M\times I$ ориентируемы, а $L_1$ и $L_2$ естественно ориентированы,
то в нормальном расслоении к $L$ есть естественная ориентация, индуцирующая естественные ориентации
$L_1$ и $L_2$.
Дополним построенную ортонормированную систему $n-1$ нормальных векторов
в каждой точке многообразия $L$ единичным вектором до ортонормированного базиса, ориентированного
положительно. В результате получим требуемое оснащение многообразия $L$.
Тем самым утверждение, выделенное курсивом, доказано.
\smallskip

{\it Завершение доказательства теоремы~\ref{d-th1} 
в случае
$w_2(M)\cdot\rho_2 H_2(M,\Z)\ne0$.}
Нам достаточно показать, что в этом случае любые два оснащенных зацепления одинаковой степени оснащенно кобордантны.
Предположим, от противного, что $L_1$ и $L_2$ --- пара (оснащенно) некобордантных зацеплений,
$L$ --- неоснащенный кобордизм между ними.
Если $\bar w_2(L)=0$, то сразу получаем противоречие (см. выделенное утверждение в начале доказательства теоремы~\ref{d-th1}
).
Поэтому предположим, что $\bar w_2(L)=1$.
Здесь $\bar w_2(L)\in H_0(L;\Z_2)\cong \Z_2$, так как $L$ связно.
Построим тогда новый кобордизм $L'$ между теми же оснащенными зацеплениями $L_1$ и $L_2$,
для которого $\bar w_2(L')=0$.
Рассмотрим элемент $\beta\in H_2(M;\Z)$, такой что $w_2(M)\cdot\rho_2\beta=1$.
Пусть $K\subset M\times\frac{1}{2}$ --- связное ориентируемое 2-подмногообразие
общего положения, реализующее класс $\beta$.
Тогда 
$$\bar w_2(K)=w_2(M)|_K=w_2(M)\cdot [K]=w_2(M)\cdot\rho_2\beta=1\mod2.$$
Действительно, первое равенство в данной цепочке равенств следует из
формулы Ву для классов Штиффеля-Уитни суммы двух расслоений:
$w_2(M)|_K=w_2(K)+w_1(K)\cdot\bar w_1(K)+\bar w_2(K)$, в которой $w_2(K)=w_1(K)=0$, так
как $K$ --- ориентируемое 2-многообразие. (Первое равенство может быть также доказано и
непосредственно.) Второе равенство следует из геометрического определения классов Штифеля-Уитни,
поскольку $K$ связно.

Положим $L'=L\sharp K$ ($L\cap K=\emptyset$ по общему положению).
Тогда $\bar w_2(L')=\bar w_2(L)+\bar w_2(K)=0$ (см. несложное утверждение~\ref{d-cl3} ниже).
Тем самым $L'$ можно оснастить. Получаем противоречие, которое доказывает теорему~\ref{d-th1}
.$\square$

\begin{claim}\label{d-cl3}
Пусть $K^2,L^2\subset M^{n+2}$ --- пара непересекающихся
связных ориентируемых подмногообразий с краем, причем на $\partial K$ and $\partial L$
заданы оснащения.
Тогда
$$
\bar w_2(K\sharp L)=\bar w_2(K)+\bar w_2(L),
$$
где группы
$H_0(X;\Z_2)$ отождествляются с $\Z_2$ для $X=K\sharp L$, $K$ и $L$.
\end{claim}


\begin{proof}[Доказательство утверждения~\ref{d-cl3}]
Рассмотрим пару малых 2-дисков
$k\subset K$ и $l\subset L$.
Пусть $kl\cong S^1\times I$ --- узкая трубка, для которой $\partial kl=\partial k\sqcup\partial l$
и $kl$ касается дисков $k$ и $l$.
Фиксируем тривиальные оснащения дисков $k$ и $l$ (и, следовательно, их границ $\partial k$ и
$\partial l$).
Из геометрического определения выше легко видеть, что
$\bar w_2(K\sharp L)=\bar w_2(K-k)+\bar w_2(kl)+\bar w_2(L-l)$.
С другой стороны, аналогично можно проверить, что
$\bar w_2(K)=\bar w_2(K-k)+\bar w_2(k)$ и
$\bar w_2(L)=\bar w_2(L-l)+\bar w_2(l)$.
Так как $\bar w_2(kl)=\bar w_2(k)=\bar w_2(l)=0$,
то из этих трех равенств получаем требуемое
$\bar w_2(K\sharp L)=\bar w_2(K)+\bar w_2(L)$.
\end{proof}


{\it Завершение доказательства теоремы~\ref{d-th1} 
в случае $w_2(M)\cdot\rho_2 H_2(M,\Z)=0$.}
Нам достаточно показать, что для фиксированного $[L_1]$ отображение $[L_2]\mapsto w_2(L)$
определено корректно и дает биекцию $\deg^{-1}\alpha\to\Z_2$.

{\it Корректность определения отображения $[L_2]\mapsto w_2(L)$.}
Пусть $L_1'$ и $L_2'$ --- пара оснащенных зацеплений в $M$,
оснащенно кобордантных $L_1$ и $L_2$ соответственно.
Пусть $L'$ --- (неоснащенный) кобордизм между ними.
Достаточно показать, что
$w_2(L)=w_2(L')$, причем только для случая, когда
$L_1$ и $L_1'$, $L_2$ и $L_2'$, $L$ и $L'$
находятся в общем положении.
Будем считать также, что $L_1,L_1'\subset M\times 1$,
$L_2,L_2'\subset M\times 0$ и $L,L'\subset M\times[0,1]$.
Изменим знак первого векторного поля оснащений зацеплений $L_1'$ и $L_2'$.
Обозначим полученные оснащенные зацепления через $-L_1'$ и
$-L_2'$ соответственно.
Обозначим через $\bar w_2(-L')$ относительный класс Штиффеля-Уитни для многообразия $L'$
c обращенными оснащениями границы $\partial L'$ (они построены выше).
Тогда $\bar w_2(-L')=-\bar w_2(L')$.
Далее, как $L_1\sqcup (-L_1')$, так и $L_2\sqcup (-L_2')$
оснащенно кобордантны нулю, то есть пустому зацеплению.
Пусть $L_{+}\subset M\times[1,+\infty)$ и
$L_{-}\subset M\times(-\infty,0]$ --- соответствующие оснащенные кобордизмы.
Тогда по определению $\bar w_2(L_{+})=\bar w_2(L_{-})=0$.
По общему положению $L\cap L'=\emptyset$.
Рассмотрим многообразие $K=L \cup L_{+}\cup L'\cup L_{-}$.
Из геометрического определения легко видеть, что
$$
\bar w_2(K)=
\bar w_2(L) + \bar w_2(L_{+}) +
\bar w_2(-L') + \bar w_2(L_{-})=
\bar w_2(L) - \bar w_2(L').
$$
Поэтому для доказательства корректности нам остается показать, что $w_2(K)=0$.
Спроектируем $K$ на $M$ с помощью 'вертикальной' проекции $M\times \Bbb R\to M$.
Пусть $\beta$ --- целочисленный гомологический класс образа многообразия $K$
при данной проекции.
Но $\bar w_2(K)=\bar w_2(M)\cdot\rho_2\beta=0$. Корректность нашего определения доказана.

{\it Инъективность отображения $[L_2]\mapsto w_2(L)$.}
Пусть $L_2'$ --- оснащенное зацепление, $L'$ --- связный (неоснащенный) кобордизм
между $L_1$ и $L_2'$.
Достаточно доказать, что если $\bar w_2(L)=\bar w_2(L')$, то $[L_2]=[L_2']$.
В самом деле, можно считать, что $L_1\subset M\times 0$,
$L_2\subset M\times 1$, $L_2'\subset M\times(-1)$,
$L\subset M\times[0,1]$ и $L'\subset M\times[-1,0]$.
Тогда $L\cup L'$ будет (неоснащенным) кобордизмом между $L_2$ и $L_2'$.
Согласно геометрическому определению выше
$\bar w_2(L\cup L')=\bar w_2(L)+\bar w_2(L')=0$,
что доказывает требуемую инъективность.

{\it Сюръективность отображения $[L_2]\mapsto w_2(L)$.}
Достаточно показать, что некоторый класс $[L_2]$ отображается в $1$.
Можно считать $L_1\cong S^1$. Фиксируем данный гомеоморфизм.
Так как $M$ ориентируемо, то существует оснащение $f_1$ кривой $L_1$.
Обозначим через $f_1(x)$ репер оснащения в точке $x\in S^1$.
Возьмем отображение $\varphi:S^1\to SO(n)$, реализующее образующую $\pi_1(SO(n))\cong\Z_2$
(что справедливо ввиду $n\ge 3$).
Определим новое оснащение $f_2$ кривой $L_1$ формулой $f_2(x)=\varphi(x)f_1(x)$.
Полученное оснащенная кривая $L_2$ --- кривая $L_1$ с оснащением $f_2$ --- искомая.
Проверим это.
Возьмем $L=L_1\times I$.
Легко видеть, что $\bar w_2(L)=1$.
Действительно, предположим противное. Тогда оснащения $f_1$ и $f_2$
продолжаются до оснащений цилиндра $L_1 \times I$.
Это оснащение задает гомотопию между $\varphi$ и постоянным отображением
$S^1\to SO(n)$, что противоречит выбору $\varphi$.
Полученное противоречие завершает доказательство.
$\square$


\bigskip

{\bf Благодарности.} Авторы благодарны А. Скопенкову и рецензенту журнала Israel Journal of Mathematics за полезные обсуждения и замечания.





\begin{thebibliography}{99}

\bibitem{CRS3} M. Cencelj, D. Repov\v{s} and M. Skopenkov,
\textit{Classification of framed links in 3--manifolds},
Proc. Indian Acad. Sci. (Math. Sci.) \textbf{117:3} (2007), p. 301--306;
arXiv:math-gt/0705.4166.

\bibitem{DNF79} Б. А. Дубровин, С. П. Новиков, А. Т. Фоменко, \textit{Современная геометрия: методы и приложения}, Наука, Москва (1979).

\bibitem{FoFu89} A. T. Фоменко и Д. Б. Фукс, \textit{Курс гомотопической топологии}, Наука, Москва 1989.

\bibitem{FU91} D.S. Freed, K.K. Uhlenbeck, \textit{Instantons and four-manifolds},
Springer-Verlag,
London, 1991.

\bibitem{Pon38} L. S. Pontryagin, \textit{Classification des transformations d'un complexe
(n+1)-dimensionel dans une sphere n-dimensionelle}, C. R. Paris \textbf{206} (1938), p. 1436--1438.

\bibitem{Pon41} L. S. Pontryagin, \textit{A classification of mappings of the 3-dimensional complex into
the 2-dimensional sphere}, Rec. Math. (Mat. Sbornik) \textbf{9:51} (1941), p. 331--363.

\bibitem{Pon76} Л. С. Понтрягин, \textit{Гладкие многообразия и их применения в теории гомотопий}, Наука, Москва 1976.

\bibitem{Ste47} N. Steenrod, \textit{Products of cocycles and extensions of mappings}, Ann. math. \textbf{48:2} (1947), p. 290--320.

\end{thebibliography}
\end{document}